\documentclass[12pt,draft]{amsart}
\usepackage{amsfonts,amscd,latexsym}
\title[On a rigidity criterion...]
{On a rigidity criterion for del Pezzo fibrations over $\POn$.}
\thanks{This work was partially supported by the grants RFBR no. 99--01--01132,
Grant of Leading Scientific Schools no. 96--15--96146, and INTAS-OPEN 97/2072}
\author{Mikhail Grinenko}
\address{Steklov Mathematical Institute}
\address{Max-Planck Institut f\"ur Mathematik}
\email{grin@mi.ras.ru / grinenko@mpim-bonn.mpg.de}
\date{}

\newtheorem{theorem}{\sc Theorem}[section]
\newtheorem{proposition}[theorem]{\sc Proposition}
\newtheorem{lemma}[theorem]{\sc Lemma}

\newtheorem{conjecture}[theorem]{\sc Conjecture}
\newtheorem{definition}[theorem]{\sc Definition}

\makeatletter

\@addtoreset{equation}{section}
\newcommand{\l@abcd}[2]{\hbox to\textwidth{#1\dotfill #2}}

\makeatother

\newcommand*{\mybegintheorem}[1]{\begin{trivlist}\it%
      \item[\hspace{\labelsep}{\bf #1}]}
\newcommand*{\myendtheorem}{\end{trivlist}}
\newenvironment*{theorem*}{\mybegintheorem{Theorem.}}{\myendtheorem}
\newenvironment*{proposition*}{\mybegintheorem{Proposition.}}{\myendtheorem}
\newenvironment*{corollary*}{\mybegintheorem{Corollary.}}{\myendtheorem}
\newenvironment*{definition*}{\mybegintheorem{Definition.}}{\myendtheorem}

\theoremstyle{remark}
\newtheorem{remark}[theorem]{\sc Remark}

\renewcommand{\phi}{\varphi}
\renewcommand{\epsilon}{\varepsilon}

\newcommand{\lra}{\longrightarrow}
\newcommand{\lla}{\longleftarrow}

\newcommand{\PQ}{{{\mathbb P}^4}}
\newcommand{\PT}{{{\mathbb P}^3}}
\newcommand{\PTw}{{{\mathbb P}^2}}
\newcommand{\POn}{{{\mathbb P}^1}}
\newcommand{\ZA}{{\mathbb Z}}
\newcommand{\QA}{{\mathbb Q}}

\newcommand{\mc}{\mathcal}

\newcommand{\eqdef}{\stackrel{\rm def}{=}}

\newcommand{\mult}{\mathop{\rm mult}\nolimits}

\newcommand{\Pic}{\mathop{\rm Pic}\nolimits}
\newcommand{\Wl}{\mathop{\rm Wl}\nolimits}

\newcommand{\chr}{\mathop{\rm char}\nolimits}

\newcommand{\eps}{\varepsilon}

\begin{document}
\abstract{We discuss the rigidity problem for Mori fibrations on del Pezzo 
surfaces of degree 1, 2 and 3 over $\POn$ and formulate the following
conjecture: a del Pezzo
fibration $V/\POn$ is rigid iff for any linear system without fixed 
components its quasi-effective and adjunction thresholds coincide.
We prove the "only if" part of this conjecture.}
\endabstract
\maketitle

Birational rigidity is one of the main notions in higher-dimensional birational 
geometry. Starting from the famous paper of V.A.Iskovskikh and Yu.I.Manin on 
three-dimensional quartics (\cite{IskMan}), in nearly all known cases 
(see \cite{Isk1}, \cite{Pukh3}) rigidity was proved using the maximal singularities 
method. Applying the same ideas to the relative situation ("over the base"), 
V.G.Sarkisov had a first success with conic bundles (\cite{Sark}) and formulated the 
essentials of "working over the base" as a program (called "Sarkisov" now). 
The case of del Pezzo fibrations turned out to be more complicated, and the first 
successful results appeared only recently in works of A.V.Pukhlikov (\cite{Pukh2}).
His technics (\cite{Pukh1}) allow to apply the maximal singularities method in any
dimension (e.g., \cite{Pukh3}). Another approach to maximal singularities is given
in \cite{Corti2}.

Not long ago, A.Corti proved the Sarkisov program in dimension 3 using Mori's 
log theory (\cite{Corti1}). An up-to-date exposition of it can be found also in
\cite{BM}. All required technical tools are given in \cite{Kollar}
(see also \cite{KM}).

In this paper, we apply the Sarkisov program to study the rigidity problem for 
3-fold Mori fibrations on del Pezzo surfaces of degree 1,2 and 3. Using all known
examples of rigid and nonrigid del Pezzo fibrations (\cite{Pukh2}, \cite{Grin1},
\cite{Grin2}), we formulate a hypothetical criterion of rigidity
(conjecture \ref{conject}) and then prove its "only if" part (theorem \ref{main_th}).

In winter 1998, during one of his lectures at the Moscow State University,
V.V.Shokurov said that it would be extremely useful if someone formulates at
least some kind of rigidity criterion for Mori fibrations. This work was 
stimulated by the suggestion of Shokurov.

This result was obtained during my stay at the Max-Planck-Institut f\"ur 
Mathematik in Bonn. I would like to use this opportunity and thank the Directors
and the staff of the Institute for hospitality. Also I would like to thank
Prof. Seunghun Lee who pointed me on some inaccuracies in the text.

\section{Preliminary.}
\label{sec1}

In this section we recall some key notions and results of the Sarkisov program.
All schemes are defined over an algebaically closed field $k$,  $\chr k=0$.

\begin{definition}
A Mori fiber space $\pi: V\to S$ is an extremal contraction of fibering type,
i.e., $V$ and $S$ are normal, $\dim S < \dim V$, $V$ has only $\QA$-factorial and
terminal singularities, $\pi_*{\mc O}_V={\mc O}_S$, $\rho(V/S)=1$ and $(-K_V)$ is
$\pi$-ample.
\end{definition}

\begin{definition}
A Mori fiber space $\pi:V\to S$ is said to be (birationally) rigid, if any
birational map $\chi:V\dasharrow V'$ to another Mori fibering $\pi':V'\to S'$ is
fiberwise, i.e., there exists a birational map $\mu:S\dasharrow S'$ such that
$\pi'\circ\chi=\mu\circ\pi$:
$$
\begin{CD}
V @. \stackrel{\chi}{\dasharrow} @. V' \\
@V{\pi}VV  @. @VV{\pi'}V \\
S @. \stackrel{\mu}{\dasharrow} @. S'
\end{CD}
$$
In this case we will often say that $V/S$ and $V'/S'$ are birational over the 
base and $\chi$ is square.
\end{definition}

In the sequel, we intend to apply the log minimal model program (LMMP) to 
log pairs of kind $K_V+\eps{\mc D}$, where ${\mc D}$ is a linear system without
fixed components. Let us formulate some statements about such pairs.

First, any Weil divisor $D$ on $V$ can be viewed as
$$
     D\sim n(-K_V)+\pi^*(A),
$$    
where $A$ is  a Weil divisor on $S$. Here $n$ is uniquely defined and said to be
{\it the quasi-effective threshold} $\mu(V,D)$ of $D$ (or a linear system 
containing $D$).

Let $D$ be a general element of ${\mc D}$. By definition, a log pair 
$K_V+a{\mc D}$, $a\in\QA_{\ge 0}$, is terminal (resp. canonical, log terminal, log 
canonical), if $K_V+aD$ is terminal (resp. canonical, log terminal, log canonical)
in the common sense (e.g., \cite{KM}). Note that if $K_V+a{\mc D}$ is canonical
then $V$ has at most terminal singularities. Moreover, running LMMP for a 
canonical pair $K+a{\mc D}$, we always stay in the class of varieties with 
$\QA$-factorial and terminal singularities (\cite{KM}, lemma 3.38, or 
\cite{Corti1}, proposition 2.6).

Let $\chi:V\dasharrow V'$ be a birational map of 3-fold Mori fibrations 
$\pi:V\to S$ and $\pi':V'\to S'$, ${\mc H'}=|n'(-K_{V'})+\pi^*(A')|$ a nonempty
linear system without base points, $A'$ an ample divisor, and
${\mc H}=\chi_*^{-1}{\mc H'}\subset|n(-K_V)+\pi^*(A)|$ the strict transform of
${\mc H'}$ on  $V$. The following facts play the crucial role in birational 
geometry  of Mori fiber spaces (\cite{Corti1}, proposition 3.5 and theorem 4.2):
\begin{proposition}
\label{Corti_th1}
(i) $n\ge n'$, and if $n=n'$ then $\chi$ is square;

\noindent (ii) if ${\mc H'}$ is ample and $K_V+\frac1n{\mc H}$ is canonical and 
nef, then $\chi$ is a biregular isomorphism over the base;

\noindent (iii) if $\chi$ is square and an isomorphism in codimension 1, then
it is biregular.
\end{proposition}

Now we outline the Sarkisov program. Each step of this program starts with
some Mori fibration $X/S$ and gives us another Mori fibration $X_1/S_1$ using
one of four types of links. Given a birational map $\chi:V\dasharrow V'$ 
between two Mori fibrations $V/S$ and $V'/S'$, we would like to decompose 
$\chi$ into a chain of such links. These links are given below (all 
contractions are extremal with respect to some log pairs; by "flips" we mean
a chain of log flips):

\begin{description}
\item[{\bf links of type 1}]
$$
\begin{array}{ccccc}
X & \gets & Z & \stackrel{flips}\dasharrow & X_1 \\
\downarrow &&&& \downarrow \\
S & \multicolumn{3}{c}{\lla} & S_1
\end{array}
$$
where $Z\to X$ is a divisorial contraction, $S_1\to S$ is a morphism
with connected fibers and $\rho(S_1/S)=1$.

\item[{\bf links of type 2}]
$$
\begin{array}{ccccccc}
X & \gets & Z & \stackrel{flips}\dasharrow & Z' & \to & X_1 \\
\downarrow &&&&&& \downarrow \\
S & \multicolumn{5}{c}{\stackrel{\sim}\lla} & S_1
\end{array}
$$
where $Z\to X$  and $Z'\to X_1$ are divisorial contractions.

\item[{\bf links of type 3}]
$$
\begin{array}{ccccc}
X & \stackrel{flips}\dasharrow & Z' & \to & X_1 \\
\downarrow &&&& \downarrow \\
S & \multicolumn{3}{c}{\lra} & S_1
\end{array}
$$
it is the inverse of type 1.

\item[{\bf links of type 4}]
$$
\begin{array}{ccccc}
X & \multicolumn{3}{c}{\stackrel{flips}\dasharrow} & X_1 \\
\downarrow &&&& \downarrow \\
S & \to & T & \gets & S_1
\end{array}
$$
where $S\to T$ and $S_1\to T$ are morphisms with connected fibers, $T$ is
normal, and $\rho(S/T)=\rho(S_1/T)=1$.
\end{description}
The decomposition procedure exists and terminates after finitely many 
steps (\cite{Corti1}, \cite{BM}).

In this paper our main subject is 3-fold Mori fibrations on del Pezzo 
surfaces of degree 1, 2 and 3 over $\POn$ (briefly "del Pezzo fibrations
of degree \ldots"). Given a del Pezzo fibration $V/\POn$, we may consider
its fiber $V_{\eta}$ over the generic point $\eta$ of $\POn$. This is a
del Pezzo surface of degree resp. 1, 2 or 3 with 
$\Pic(V_{\eta})\simeq\ZA[-K_{V_{\eta}}]$. Viewing $V_{\eta}$ as a 
two-dimensional Mori fibration, we may apply the Sarkisov program. 
The result is copmletely described in \cite{Isk3}, theorem 2.6:
\begin{proposition}
\label{Isk_prop}
Let $X/S$ and $X_1/S_1$ be two-dimensional Mori fibrations (over a perfect
nonclosed field), $X$ a del Pezzo surface of degree $d=(-K_X)^2$ 
(so $S$ is a point). Then:

\noindent i) there are no links of types 1 and 3 between $X/S$ and $X_1/S_1$;

\noindent ii) there exist links of type 2; in this case $X\simeq X_1$
 always, $d=2$ or $3$, and a link is a Geiser or Bertini involution on $X$.
\end{proposition}

Now consider a del Pezzo fibration $\pi:V\to\POn$ of degree 1, 2 or 3. 
Then the Weil divisors classes group is
$$
     \Wl(V)=\ZA[-K_V]\oplus\ZA[F],
$$
where $F$ is a class of fiber, $F=\pi^*(pt.)$. Given an effective divisor
$D$, we may define {\it the adjunction threshold} $\alpha(V,D)$. 
By definition, this is a (rational) number such that $|aD+bK_V|=\emptyset$ 
for all positive integer $a$, $b$ such that $\frac{b}{a}>\alpha(V,D)$.
It is obvious that
for a nonempty linear system ${\mc D}\subset|n(-K_V)+mF|$ we have 
$\alpha(V,{\mc D})\le n=\mu(V,{\mc D})$, and the equality holds 
only if $m\ge 0$.

I propose the following criterion of rigidity for del Pezzo fibrations:
\begin{conjecture}
\label{conject}
Let $V/\POn$ be a 3-fold Mori fibration on del Pezzo surfaces of degree 1, 2
or 3. $V/\POn$ is birationally rigid iff for any nonempty linear system 
${\mc D}$ without fixed components $\alpha(V,{\mc D})=\mu(V,{\mc D})$ 
(in other words, ${\mc D}\subset|n(-K_V)+mF|$ for $m\ge 0$).
\end{conjecture}

\begin{remark}
This conjecture is completely proved for smooth del Pezzo fibrations 
of degree 1 (\cite{Grin2}). Moreover, any smooth "general" (i.e., 
sufficiently twisted over the base) fibration of degree 2 and 3 
is rigid and satisfies the necessary condition
while it fails for all known nonrigid cases (\cite{Pukh2}, \cite{Grin2}).
\end{remark}

\begin{remark}
It does not make sense to speak about rigidity for del Pezzo fibrations of 
degree more than 3. Indeed, all such fibrations of degree 5 or higher are
rational, as it follows from \cite{Manin}. In case of degree 4 we always
 have a structure of conic bundle (\cite{Isk3}, theorem 2.6).
\end{remark}

\section{The main result.}

In this section we prove the "only if" part of conjecture \ref{conject}.

\begin{theorem}
\label{main_th}
Let $V/\POn$ be a 3-fold Mori fibration on del Pezzo surfaces of 
degree 1, 2 or 3. If $V/\POn$ is rigid, then any linear system ${\mc D}$ 
without fixed components is a subsystem of $|n(-K_V)+mF|$ for some $m\ge 0$.
\end{theorem}

We first prove an easy lemma about del Pezzo fibrations of degree 1 and 2.
Recall that there exists a morphism of degree 2 from a smooth del Pezzo surface
of degree $d=1$ or $2$ onto a plane (in case $d=2$) or a quadratic cone in $\PT$
($d=1$). This morphism is determined by the linear system $|-K|$ ($d=2$) or $|-2K|$
($d=1$). We can transpose the sheets of this cover and thus obtain a biregular
involution of such a surface.
\begin{lemma}
\label{conj_lem}
Let $V/\POn$ be a del Pezzo fibration of degree 1 or 2, $V_{\eta}$ the
generic fiber and $\tau_{\eta}$ its involution described above. Then there
exists a biregular involution of $V$ that is an extension of $\tau_{\eta}$.
\end{lemma}
\noindent{\sc Proof.} Let $\tau'$ be any extension of $\tau_{\eta}$ to the
entire $V$, so $\tau'$ is a fiberwise birational automorphism and biregular 
over an open subset of $\POn$. We have the commutative diagram
$$
\begin{array}{ccc}
V  & \stackrel{\tau'}\dasharrow & V \\
\downarrow & & \downarrow \\
\POn & \stackrel{\sim}\lra & \POn
\end{array}
$$
Let ${\mc H}=|n(-K_V)+mF|$ be a very ample linear system, $m>0$, and 
${\mc D}=\tau'{}_*^{-1}{\mc H}$ its strict transform. As it follows from
\cite{Corti1}, theorem 5.4, after finitely many links of type 2 (links of type 1
are impossible because of proposition \ref{Isk_prop}) we obtain a (square) 
birational map
$$
\begin{array}{ccc}
U  & \stackrel{\psi}\dasharrow & V \\
\downarrow & & \downarrow \\
\POn & \stackrel{\sim}\lra & \POn
\end{array}
$$
and a linear system ${\mc D}_U=\psi_*^{-1}{\mc D}$ such that $K_U+\frac1n{\mc D}_U$
is canonical. Since $\tau'\circ\psi:U/\POn\dasharrow V/\POn$ is an isomorphism
over an open subset of $\POn$, we have ${\mc D}_U\subset|n(-K_U)+lF|$. It remains
to show that $l\ge 0$. Indeed, take a common resolution 
$$
\begin{array}{ccccc}
&& W && \\
& \stackrel{f}{\swarrow} && \stackrel{g}{\searrow} & \\
U & \multicolumn{3}{c}{\stackrel{\tau'\circ\psi}{\dasharrow}}& V
\end{array}
$$
Denoting ${\mc H}_W$ the strict transform of ${\mc H}$, we get
$$
K_W+\frac1n{\mc H}_W=g^*(\frac{m}nF)+\sum a_iE_i=f^*(\frac{l}nF)+\sum b_jQ_j,
$$
where $E_i$ and $Q_j$ are exeptional divisors, all $a_i$ and $b_j$ are nonnegative.
Thus we see that $l\ge 0$. 

Now, since $K_U+\frac1n{\mc D}_U$ is canonical and nef, $\tau'\circ\psi$ is an 
isomorphism (proposition \ref{Corti_th1}) that we need. Lemma is proved.

\medskip

\noindent{\sc Proof of theorem \ref{main_th}}. Let $\pi:V\to\POn$ be a del Pezzo
fibration of degree 1, 2 or 3. Assume the converse, i.e., for $n\gg 0$ there 
exist linear systems $|n(-K_V)-F|$ without fixed components.

We claim that the log pair $K_V+\frac1n{\mc D}$ is canonical along any 
horizontal (i.e., covering the base) curve, where ${\mc D}=|n(-K_V)-F|$.
Indeed, it is always true in the case of degree 1, as it follows from 
\cite{Pukh2}, \S 3, or proposition \ref{Isk_prop}. For degree 2, if
$K_V+\frac1n{\mc D}$ is not canonical along $C$, then $C$ is a section and
there exists a section $C'$ that is conjugated to $C$ with respect to 
the biregular involution stated in \ref{conj_lem} (see \cite{Pukh2}, \S 3).
So $K_V+\frac1n{\mc D}$ is not canonical along $C'$ too, which is impossible.

Let us consider the case of degree 3. Suppose that a general member $D\in{\mc D}$ 
has multiplicity $\nu>n$ along $C$. Following to \cite{Pukh2}, we see that 
$C$ is either a section or a bi-section. We will obtain a contradiction by 
presenting an element $D'\in{\mc D}$ with $\mult_C D'\le n$.

Let $C$ be a section, $\nu=\mult_C D>n$. There exists the unique irreducible 
divisor $G\in|-2K_V+bF|$ with multiplicity 2 along $C$, where $b$ is 
uniquely determined. We apply the Geiser involution $\tau_C$ with respect to 
$C$. Blowing up $\phi:\tilde V\to V$ the curve $C$, we may consider $\tau_C$ as
an isomorphism in codimension 1 on $\tilde V$. We observe that $\tau_C$ 
transposes the strict transorm of $G$ and the exeptional divisor of $\phi$.
Taking into account proposition \ref{Isk_prop}, (ii), this allows us to compute
$\tau_C$-action on divisors. It is easy to check that
$$
   T\eqdef (\tau_C)_*^{-1}(D)\in|(2n-\nu)(-K_V)+(b(n-\nu)-1)F-(3n-2\nu)C|,
$$
i.e., $\mult_C T=3n-2\nu$. Put $D'=T+(\nu-n)G$. We get $D'\in{\mc D}$ and
$\mult_C D'=n$, a contradiction.

The same reason works if $C$ is a bi-section. There is the unique ireducible
divisor $G\in|-4K_V+bF|$ with multiplicity 5 along $C$, $b$ is uniquely 
determined. We have
$$
 T\eqdef (\tau_C)_*^{-1}(D)\in|(5n-4\nu)(-K_V)+(b(n-\nu)-1)F-(6n-5\nu)C|,
$$
and $D'=T+(\nu-n)G$ is in ${\mc D}$ and has multiplicity $n$ along $C$.

So we may assume that ${\mc D}=|n(-K_V)-F|$ has no base components and 
$K_V+\frac1n{\mc D}$ is canonical along horizontal curves.

Following the Sarkisov program, we may produice a chain of links of type 2
$$
\begin{array}{ccccccc}
V & \dasharrow & V_1 & \dasharrow & \ldots & \dasharrow & V_k \\
\downarrow & & \downarrow & & & & \downarrow \\
\POn & \lla & \POn & \lla & \ldots & \lla & \POn
\end{array}
$$
such that $K_{V_k}+\frac1n{\mc D}_k$ is canonical, where ${\mc D}_k$ is the 
strict transform of ${\mc D}$ on $V_k$. Let us note that any link in the chain
is constructed as follows. First we blow up something in a fiber, then make some
log flips over the base, and contract the strict transform of this fiber. Thus
this chain gives us an isomorphism over an open subset of $\POn$. It is clear
that ${\mc D}_k\subset|n(-K_{V_k})-aF|$, and all we need is to show that $a>0$. 
But we know that Iitaka's dimension of log divisors is a birational invariant.
Obviously, $k(V, K_V+\frac1n{\mc D})=-\infty$, and if $a\le 0$ then
$k(V_k, K_{V_k}+\frac1n{\mc D}_k)\ge 0$ since $K_{V_k}+\frac1n{\mc D}_k$ 
is canonical.
Now we can apply a link of type 3 or 4 to $(V_k, K_{V_k}+\frac1n{\mc D}_k)$.
We get either $\QA$-Fano (link of type 3), or conic bundle or del Pezzo fibration
(link of type 4). Any of these Mori structures is different from $V/\POn$.
This completes the proof of theorem \ref{main_th}.

\section{Remarks.}

\noindent{\bf 1. Projective models.} Let $\pi:V\to\POn$ be a Gorenstein del
Pezzo fibration of degree $d=1$, $2$ or $3$. Then we can present a very 
convenient projective model of $V/\POn$ constructed as follows (see \cite{Isk2}).

First, we must deal with an arbitrary Gorenstein del Pezzo surface.
\begin{proposition}
Let $S$ be a reduced irreducible Gorenstein del Pezzo surface of degree
$d=K_S^2$. Then
we have the following projective models:
\begin{description}
\item[$d=1$] the linear system $|-2K_S|$ is base point free and gives us a
morphism $\phi:S\to Q$ to a quadratic cone $Q\subset\PT$ branched over a cubic
section that does not pass through the cone vertex. The unique base point of
$|-K_S|$ lies over this vertex;
\item[$d=2$] the linear system $|-K_S|$ is base point free and gives us a
morphism $\phi:S\to\PTw$ branched over a quartic;
\item[$d=3$] the linear system $|-K_S|$ is base point free and embeds $S$ into
$\PT$.
\end{description}
\end{proposition}
\noindent{\sc Proof.} The normal case is proved in \cite{HiWa}, corollary 4.5 
and proposition 4.6. Non-normal surfaces are completely described in \cite{Reid},
section 1.4 and corollary 4.10, or \cite{Fu}.

Suppose $d=1$. For $m\gg 0$ the linear system $|-2K_S+mF|$ gives a morphism 
$\phi:V\to Q$ of degree 2 onto a variety $Q$ fibered on quadratic cones over $\POn$.
We may assume 
$$Q\subset P\eqdef
{\bf Proj}_{\POn}({\mc O}\oplus{\mc O}(n_1)\oplus{\mc O}(n_2)\oplus{\mc O}(n_3)),
$$
where $n_3\ge n_2\ge n_1\ge 0$.
Let $M$ be the class of tautological bundle on $P$, $L$ the class of fiber,
$t_0$ a section in $P$ such that $t_0\circ M=0$, and $l$ the class of a line in 
fiber. Suppose also $t_b\sim t_0+\eps l$ is the curve of cone vertices in $Q$, 
$R$ is the ramification divisor of $\phi$, 
$G\sim M^2-(n_2+n_3)ML$ is the class of the "minimally twisted" surface in $Q$,
$G_V=\phi^*(G)$, and $H=\phi^*(M)$. The following fact is proved in \cite{Grin2}, 
lemmas 2.1--2.3 (it does work in singular cases too):
\begin{proposition} The only following cases may occur: either

\noindent (i) $\eps=0$, and then $2n_2=n_1+n_3$, $Q\sim 2M-2n_2L$, $R\sim 3M$,
$-K_V\sim G_V+(2-\frac12n_1)F=\frac12H+(2-n_2)F$,

or

\noindent (ii) $\eps=n_1>0$, and then $n_3=2n_2$, $n_2\ge 3n_1$, 
$Q\sim 2M-2n_2L$, $R\sim 3M-3n_1L$, 
$-K_V\sim G_V+(2+\frac12n_1)F=\frac12H+(2+\frac12n_1-n_2)F$.
\end{proposition}

Now let $d=2$. The linear system $|-K_V+mF|$ gives the degree 2 morphism
$$
     \phi: V\to {\bf Proj}_{\POn}({\mc O}\oplus{\mc O}(n_1)\oplus{\mc O}(n_2)),
$$
where $n_2\ge n_1\ge 0$, with the ramification divisor $R$. Suppose $b=n_1+n_2$.
Using the same notation as before, we can see that $-K_V\sim H+(2-a-b)F$, where
$a$ is determined by $R\sim 4M+aL$ (\cite{Isk2}, section 1.3, or \cite{Grin2}, 
section 3.1).

Finally, the case of degree $d=3$ can be viewed in the same way 
(\cite{Isk2}, section 1.3).

\noindent{\bf 2. Proof of \ref{main_th} in Gorenstein cases.}
Let us consider a Gorenstein del Pezzo fibration $V/\POn$
of degree 1 or 2. Using the projective models described above, we may easily 
observe that
\begin{center}
{\it if there exists a nonempty linear system $|n(-K_V)-F|$ without fixed 
components, so does $|-K_V-F|$ (with no base divisors).}
\end{center}
Moreover, a general element of $|-K_V-F|$ is a normal elliptic fibration
over $\POn$ without multiple fibers, so it is rational (use the adjunction 
formula). Pick a pencil ${\mc P}\subset|-K_V-F|$ and take a resolution
$\psi:\tilde V\to V$ of its base points. Then ${\mc P}$ defines a morphism
$\tilde V\to\POn$, and it only remains to apply the relative minimal model 
program. We get another structure of Mori fibration (either del Pezzo or 
conic bundle).

Note that it fails in the case $d=3$. Here is a counterexample. Let 
$X\subset\PQ$ be a general quartic containing a plane $P$. Blowing up $P$ 
(as a surface in $\PQ$), we obtain a small resolution $\psi:V\to X$ of 
singularities of $X$. Obviously, $V$ is a del Pezzo fibration of degree 3.
Let $E$ be the strict transform of $P$. The linear system $|-3K_V-F|$ has no
base divisors, but $E$ is the unique element of $|-K_V-F|$.
By the way,  as it follows from \ref{main_th}, another 
structure of Mori fibration must exist.
We can find it easily. Note that $\psi|_E:E\to P$ contracts 9 lines
to 9 points of some elliptic pencil on $P$. Making a flop at these lines 
simultaneously and then contracting the strict transform of $E$ (it is a 
plane with the normal bundle ${\mc O}(-2)$), we obtain a $\QA$-Fano 
variety with a singular point of index 2.

\end{document}